\documentclass[11pt]{article}
\usepackage{amsfonts,latexsym,amsmath,amscd,geometry}
\usepackage{amssymb}
\usepackage{latexsym}

\newcommand \nc{\newcommand}
\newtheorem{theorem}{Theorem}[section]
\newtheorem{lemma}[theorem]{Lemma}

\newtheorem{corollary}[theorem]{Corollary}
\newtheorem{definition}[theorem]{Definition}

\newtheorem{remark}[theorem]{Remark}

\nc{\ba}{\begin{array}}\nc{\ea}{\end{array}}
\nc{\be}{\begin{eqnarray}}\nc{\ee}{\end{eqnarray}}
\nc{\beq}{\begin{equation}}\nc{\eeq}{\end{equation}}
\nc{\bex}{\begin{eqnarray*}}\nc{\eex}{\end{eqnarray*}}
\nc{\btm}{\begin{theorem}} \nc{\etm}{\end{theorem}}
\nc{\blm}{\begin{lemma}} \nc{\elm}{\end{lemma}}
\nc{\R}{\mathbb{R}} \nc{\va}{\varepsilon} \nc{\ls}{\limits}

\def\de{\Delta}\def\pa{\partial}\def\om{\Omega}
\def\pf{\noindent{\bf Proof.\quad}}\def\endpf{\hfill$\Box$}
\def\les{\lesssim}\def\di{\mbox{div\,}}
\def\curl{\mbox{curl\,}}
\newcommand \qed {\hfill $\Box$}

\begin{document}
\title{Blow up criterion for compressible nematic liquid crystal flows in dimension three}
\author{Tao Huang\footnote{Department of Mathematics, University of Kentucky,
Lexington, KY 40506, USA.} \quad Changyou Wang$^*$ \quad Huanyao Wen\footnote{School of Mathematical Sciences, South China Normal University, Guangzhou, 510631, P. R. China.}}
\maketitle

\begin{abstract} In this paper, we consider the short time strong solution to a simplified hydrodynamic flow modeling the compressible, nematic liquid crystal materials in dimension three. We stablish a criterion
for possible breakdown of such solutions at finite time in terms of the temporal integral of both the maximum norm of the deformation tensor of velocity gradient and the square of maximum norm of gradient of liquid crystal director field.
\end{abstract}

\section {Introduction}
\setcounter{equation}{0}
\setcounter{theorem}{0}

Nematic liquid crystals are aggregates of molecules which possess same orientational order and are made of elongated,
rod-like molecules. The continuum theory of liquid crystals was developed by Ericksen \cite{Ericksen}  and Leslie \cite{Leslie}
during the period of 1958 through 1968, see also the book by de Gennes \cite{Gennes}.  Since then there have been remarkable
research developments in liquid crystals from both theoretical and applied aspects.  When the fluid containing nematic liquid crystal materials is at rest, we have the well-known Ossen-Frank theory for static nematic liquid crystals. The readers can refer to the poineering work by Hardt-Lin-Kinderlehrer \cite{HLK} on the analysis of energy minimal configurations of namatic liquid crystals.  In general, the motion of fluid always takes place. The so-called Ericksen-Leslie system is a macroscopic continuum description of the time evolution of the materials under the influence of both the flow velocity field $u$ and the macroscopic description of the microscopic orientation configurations $d$ of rod-like liquid crystals. 

When the fluid is an incompressible, viscous  fluid, Lin \cite{Lin} first derived a simplified Ericksen-Leslie equation modeling liquid crystal flows
in 1989. Subsequently, Lin and Liu \cite{LL, LL1} made some important analytic studies, such as the existence of weak and strong solutions and the partial regularity of suitable solutions,  of the simplified Ericksen-Leslie system,
under the assumption that the liquid crystal director field is of varying length by Leslie's terminology or variable degree of orientation by Ericksen's terminology.  

When the fluid is allowed to be compressible, the Ericksen-Leslie system becomes more complicate and there seems very few analytic works available yet.  We would like to mention that very recently, there have been both modeling study, see Morro \cite {Morro},  and numerical study, see Zakharov-Vakulenko \cite{ZV}, on the hydrodynamics of compressible nematic liquid crystals  under the influence of temperature gradient or electromagnetic forces.   

The main aim of this paper,  and the companion paper \cite{huang-wang-wen} as well,
is an attemption to initiate some analytic study for the flow of compressible nematic liquid crystals. We will
mainly address several issues on the strong solutions.  More precisely, we will focus on the blow-up criterion on
strong solutions in this paper.

Now we start to describe the problem. 
Let $\Omega\subset\mathbb R^3$ be either a bounded smooth domain or the entire $\mathbb R^3$, 
we will consider a simplified version of Ericksen-Leslie equation that models the hydrodynamic flow of
compressible, nematic liquid crystals in $\Omega$:
\begin{align}
\rho_t+\nabla \cdot (\rho u)&=0, \label{clc-1} \\
\rho u_t+\rho u\cdot\nabla u+\nabla (P(\rho))&=\mathcal L u-\de d\cdot\nabla d,\label{clc-2}\\
d_t+u\cdot\nabla d&=\de d+|\nabla d|^2d,\label{clc-3}
\end{align}
where $\rho:\om\to\mathbb R_+$ is the density of the fluid,\ $u:\om\to\mathbb R^3$ is the fluid velocity field, $P(\rho):\Omega\to\mathbb R_+$ denotes the pressure of the fluid,  $d:\Omega\to S^2$
represents  the macroscopic average of the nematic liquid crystal orientation field, $\nabla\cdot (={\rm{div}})$
denotes the divergence operator on $\mathbb R^3$, and $\mathcal L$ denotes the Lam$\acute{\rm e}$ operator defined by
\begin{equation}\label{lame}
\mathcal Lu=\mu\de u+(\mu+\lambda)\nabla(\nabla\cdot u),
\end{equation}
where $\mu$ and $\lambda$ are the shear viscosity and the bulk viscosity coefficients
of the fluid repsectively,  which are assumed to satisfy the following physical condition:
\beq\label{viscosity}
\mu>0, \quad 3\lambda+2\mu\geq 0.
\eeq
The pressure $P(\rho)$, as a given continuous function of $\rho$, is usually determined by the equation of states. Through this paper, we assume that
\beq\label{regularity_p}
P:[0,+\infty)\to\mathbb R\ {\rm{is\ a\  locally\  Lipschitz\ continuous\ function}}.
\eeq
Notice that (\ref{clc-1}) is the equation for conservation of mass, (\ref{clc-2}) is the linear momemtum equation,
and (\ref{clc-3}) is the angular momentum equation.  We would like to point out that the system (\ref{clc-1})-(\ref{clc-3})
includes several important equations as special cases:

(i) When $\rho$ is constant, the equation (\ref{clc-1}) reduces
to the incompressibility condition of the fluid ($\nabla\cdot u=0$), and the system (\ref{clc-1})-(\ref{clc-3}) becomes
the equation of incompressible flow of namatic liquid crystals provided that $P$ is a unknown pressure function. This
was previously proposed by Lin \cite{Lin} as a simplified Ericksen-Leslie equation modeling incompressible liquid crystal flows.

(ii) When $d$ is a constant vector field,  the system (\ref{clc-1})-(\ref{clc-2})
becomes a compressible Navier-Stokes equation, which is an extremely important equation to describe 
motion of compressible fluids.
It has attracted great interests among many analysts and there have been many important developments (see,
for example,  Lions \cite{Lions2}, Feireisl \cite{Feireisl} and references therein).

(iii) When both $\rho$ and $d$ are constants, the system (\ref{clc-1})-(\ref{clc-2}) becomes the incompressible Naiver-Stokes
equation provided that $P$ is a unknown pressure function, the fundamental equation to describe Newtonian fluids (see, Lions \cite{Lions1}
and Temam \cite{Temam} for survey of important developments).

(iv) When $\rho$ is constant and $u=0$, the system (\ref{clc-1})-(\ref{clc-3}) reduces to the equation for heat flow of harmonic maps into $S^2$.  There have been extensive studies on the heat flow of harmonic maps in the past few decades (see, for example,
the monograph by Lin-Wang \cite{LW1} and references therein).

From the viewpoint of partial differential equations, the system (\ref{clc-1})-(\ref{clc-3}) is a highly nonlinear system coupling between hyperbolic equations and parabolic equations.
It is very challenging to understand and analyze such a system, especially when the density function $\rho$ may vanish
or the fluid takes vacuum states.

In this paper,  the system (\ref{clc-1})-(\ref{clc-3}) will be studied along with the initial condition:
\begin{align}
(\rho, u, d)\Big|_{t=0}&= (\rho_0, u_0, d_0), \label{clcinitial}
\end{align}
and one of the following three types of boundary conditions:\\
(i) Cauchy problem:
\begin{align}
\om=\R^3, \ \rho, u \mbox{ vanish and } d \mbox{ is constant at infinity (in some weak sense)}.\label{clcboundary1}
\end{align}
(ii) Dirichlet and Neumann boundary condition for $(u,d)$: $\om\subset\R^3$ is a bounded smooth domain, and
\begin{align}
(u,\ \frac{\partial d}{\partial\nu})\Big|_{\pa\Omega}&= 0, \label{clcboundary2}
\end{align}
where $\nu$ is the unit outward normal vector of $\pa\om$.\\
(iii) Navier-slip and Neumann boundary condition for $(u,d)$: $\om\subset\R^3$ is a bounded, simply connected, smooth domain, and
\begin{align}
(u\cdot\nu,\ \curl u\times\nu,\ \frac{\partial d}{\partial\nu})\Big|_{\pa\Omega}&= 0,\label{clcboundary3}
\end{align}
where $\curl u=\nabla\times u$ denotes the vorticity field of the fluid.

In order to state the definition of strong solutions to the initial and boundary value
problem  (\ref{clc-1})-(\ref{clc-3}), (\ref{clcinitial}) together with (\ref{clcboundary1})
or (\ref{clcboundary2}) or (\ref{clcboundary3}), we introduce some notations. 

We denote 
$$\int f\,dx =\int_\Omega f \,dx.$$
For $1\le r\le \infty$, denote the $L^r$ spaces and the standard Sobolev spaces as follows:
$$L^r=L^r(\Omega),  \ D^{k,r}=\left\{ u\in L^1_{\rm{loc}}(\Omega): \|\nabla^k u \|_{L^r}<\infty\right\},$$
$$W^{k,r}=L^r\cap D^{k,r},  \ H^k=W^{k,2}, \ D^k=D^{k,2},$$
$$D_0^1=\Big\{u\in L^6: \ \|\nabla u\|_{L^2}<\infty,\
{\rm{ and \ satisfies}}\ (\ref{clcboundary1}) \ {\rm{or}}\ (\ref{clcboundary2}) \ {\rm{or}}\ (\ref{clcboundary3}) 
\ {\rm{ for\ the\ part\ of }}\ u\Big\},$$
$$H_0^1=L^2\cap D_0^1, \ \|u\|_{D^{k,r}}=\|\nabla^k u\|_{L^r}.$$
Denote $$Q_T=\Omega\times [0,T] \ (T>0),$$
and let
$$\mathcal D(u)=\frac12\left(\nabla u+(\nabla u)^t\right)$$
denote the deformation tensor, which is the symmetric part of the velocity gradient.

\begin{definition} For $T>0$, $(\rho, u,d)$ is called a strong solution to the compressible nematic liquid crystal flow (\ref{clc-1})-(\ref{clc-3}) in $\Omega\times (0,T]$, if for some $q\in (3, 6]$, 
\bex &0\le \rho\in
C([0,T];W^{1,q}\bigcap H^1),\ \rho_t\in C([0,T];L^2\bigcap L^q);&\\
& u\in C([0,T];D^2\bigcap D^1_0)\bigcap L^2(0,T;D^{2,q}),\ u_t\in
L^2(0,T;D^1_0),\
 \sqrt{\rho} u_t\in
L^\infty(0,T;L^2);&\\
& \nabla d\in C([0,T];H^2)\bigcap L^2(0,T;
H^3),\ d_t\in C([0,T];H^1)\bigcap L^2(0,T; H^2),\
 |d|=1 \ {\rm{in}}\ \overline{Q}_T;& \eex
and $(\rho,u, d)$ satisfies (\ref{clc-1})-(\ref{clc-3}) a.e. in $\Omega\times (0,T]$.

\end{definition} 

For the existence of local strong solutions associated with the three types of boundary conditions, we have obtained
the following theorem in the paper \cite{huang-wang-wen}.

\begin{theorem}\label{App-thm:local}
Assume that $P$ satisfies (\ref{regularity_p}), $\rho_0\geq0$, $\rho_0\in W^{1,q}\bigcap H^1\bigcap L^1$ for some
$q\in(3,6]$, $u_0\in D^2\bigcap D^1_0$, $\nabla d_0\in H^2$ and
$|d_0|=1$ in $\overline{\Omega}$. If, in additions, the following compatibility
condition 
\be\label{first3.1} \mathcal Lu_0-\nabla
(P(\rho_0))-\de d_0\cdot\nabla d_0=\sqrt{\rho_0}g\
{\rm{ for\ some}}\ g\in L^2(\Omega,\mathbb R^3)
\ee
holds,  then there exist a positive time $T_0>0$
 and a unique strong solution $(\rho,u, d) $ of (\ref{clc-1})-(\ref{clc-3}),
(\ref{clcinitial}) together with (\ref{clcboundary1}) or 
(\ref{clcboundary2}) or (\ref{clcboundary3}) in $\Omega\times (0,T_0]$.

\end{theorem}

We would like to point out that  an analogous existence theorem of local strong solutions to the isentropic compressible Naiver-Stokes equation, under the first two boundary conditions (\ref{clcboundary1}) and
(\ref{clcboundary2}), has been previously established  by Choe-Kim \cite{CK} and  Cho-Choe-Kim \cite{CCK}. A byproduct of our theorem \ref{App-thm:local} also yields the existence of local strong solutions to the isentropic compressible Navier-Stokes
equation under the Navier-slip boundary condition (\ref{clcboundary3}).  

In dimension one, Ding-Lin-Wang-Wen \cite{DLWW} have proven that the local strong solution to
(\ref{clc-1})-(\ref{clc-3}) under (\ref{clcinitial}) and (\ref{clcboundary2}) is global.  For dimensions at least two, it is reasonable to believe that the local strong solution to (\ref{clc-1})-(\ref{clc-3}) may cease to exist globally.
In fact, there exist finite time singularities of the (transported) heat flow of harmonic maps (\ref{clc-3}) in
dimensions two or higher (we refer the interested readers to \cite{LW1} for the exact references). An important question to ask would be what is the main mechanism of possible break down of local strong (or smooth)
solutions. 

Such a question has been studied for the incompressible Euler equation or the Navier-Stokes equation by Beale-Kato-Majda in their poineering work \cite{BKM}, which showed that the $L^\infty$-bound of vorticity $\nabla\times u$ must blow up. Later, Ponce \cite{Ponce} rephrased the BKM-criterion in terms of the deformation tensor $\mathcal D(u)$.  

When dealing with the {\it isentropic}\footnote{namely, $P(\rho)=a\rho^\gamma$ for some $a>0$ and
$\gamma>1$.} compressible Navier-Stokes equation, there have recently been several very interesting works on the blow up criterion. For example,
if $0<T_*<+\infty$ is the maximum time for strong solution, then (i) Huang-Li-Xin \cite{HLX0}  established
a Serrin type criterion: 
$\lim_{T\uparrow T_*} \big(\|{\rm{div}} u\|_{L^1(0,T; L^\infty)}+\|\rho^\frac12 u\|_{L^s(0,T; L^r)}\big)=\infty$ for $\frac{2}{s}+\frac{3}{r}\le 1, \ 3<r\le\infty$; (ii) Sun-Wang-Zhang \cite{Sun-Wang-Zhang}, and independently \cite{HLX0}, showed that if $7\mu>\lambda$, then 
$\lim_{T\uparrow T_*} \|\rho\|_{L^\infty(0,T; L^\infty)}=\infty$; and (iii) Huang-Li-Xin \cite{HLX} showed
that $\lim_{T\uparrow T_*} \|\mathcal D(u)\|_{L^1(0,T; L^\infty)}=\infty$.  
 
When dealing the heat flow of harmonic maps (\ref{clc-3}) (with $u=0$),  Wang \cite{Wang} obtained
a Serrin type regularity theorem, which implies that if $0<T_*<+\infty$ is the first singular time for
local smooth solutions, then $\lim_{T\uparrow T_*} \|\nabla d\|_{L^2(0,T; L^\infty)}=\infty$.

When dealing with the incompressible nematic liquid crystal flow,  Lin-Lin-Wang \cite{Lin-Lin-Wang}  and Lin-Wang \cite{LW} have established the global existence of a unique "almost strong" solution\footnote{that has at most finitely many possible singular time.} for the initial-boundary value problem in  bounded domains
in dimension two, see also Hong \cite{Hong} and Xu-Zhang \cite{Xu-Zhang} for some related works.   
In dimension three,  for the incompressible nematic liquid crystal flow Huang-Wang \cite{HW} have obtained a BKM type blow-up criterion very recently,  while the existence of global weak solutions still
remains to be a largely open question. 

Motivated by these works on the blow up criterion of local strong solutions to the Navier-Stokes equation and
the incompressible nematic liquid crystal flow, we will establish in this paper the following blow-up criterion of breakdown of local strong solutions in finite time. 

\btm\label{umaintheorem}{\it 
Let $(\rho, u, d)$ be a strong solution of the initial boundary problem
(\ref{clc-1})-(\ref{clc-3}), (\ref{clcinitial}) together with (\ref{clcboundary1})
or (\ref{clcboundary2}) or (\ref{clcboundary3}). Assume that $P$ satisfies (\ref{regularity_p}), and the initial data $(\rho_0, u_0,d_0)$
satisfies (\ref{first3.1}). If $0<T_*<+\infty$ is the maximum time of existence,
then  
\beq\label{clcblpcondition}
\int_0^{T_*} \left(\|\mathcal{D}(u)\|_{L^{\infty}}+\|\nabla
d\|^2_{L^{\infty}}\right)\,dt=\infty. \eeq}\etm
We would like to make a few comments now.
\begin{remark} {\rm{(a) In \cite{huang-wang-wen}, we also obtained a blow-up criterion of (\ref{clc-1})-(\ref{clc-3}) under the initial condition (\ref{clcinitial}) and (\ref{clcboundary1})
or (\ref{clcboundary2}) in terms of $\rho$ and $\nabla d$. Namely, if $7\mu>9\lambda$ and
$0<T_*<+\infty$ is the maximum time of existence of strong solutions, then
$$\lim_{T\uparrow T_*} \Big(\|\rho\|_{L^\infty(0,T; L^\infty)}+\|\nabla d\|_{L^3(0,T; L^\infty)}\Big)=+\infty.$$
(b) For compressible liquid crystal flows without the nematicity constraint ($|d|=1$)\footnote{the right hand side of equation (\ref{clc-3}) is replaced by $\Delta d+f(d)$ for some smooth function $f:\mathbb R^3
\to\mathbb R^3$, e.g. $f(d)=(|d|^2-1)d$.}, Liu-Liu \cite{Liu-Liu}
have recently obtained a Serrin type criterion on the blow-up of strong solutions.\\
(c) It is a very interesting question to ask whether there exists a global weak solution to the initial-boundary value problem of (\ref{clc-1})-(\ref{clc-3}) in dimensions at least two. In dimension one, such an
existence has been obtained by Ding-Wang-Wen \cite{DWW}.}}
\end{remark}

We conclude this section by introducing the main ideas of the proof, some of which are inspired by
some of the arguments on the isentropic compressible Navier-Stokes equation by \cite{HLX}
and \cite{Sun-Wang-Zhang}.\\ (1)  It is well-known that
the bound of $\displaystyle\|\mathcal D(u)\|_{L^1_tL^\infty_x}$ yields that 
$\displaystyle\|\rho\|_{L^\infty_t L^\infty_x}$ is bounded from the equation  (\ref{clc-1}). See Lemma 2.1. \\
(2) We observe that in the equation (\ref{clc-3}), the bound of $\displaystyle(\|\mathcal D(u)\|_{L^1_tL^\infty_x}+\displaystyle\|\nabla d\|_{L^2_tL^\infty_x})$ yields that $\displaystyle\|\nabla d\|_{L^\infty_tL^r_x}$
is bounded for any $2\le r<+\infty$, which is a crucial ingredient to obtain higher order estimates
of $\rho, u, d$. See Lemma 2.3.\\ 
(3) Due to the possible vacuum state of $\rho$, the strong nonlinearities of the convection
term $u\cdot\nabla u$ and the induced stress tensor $\Delta d\cdot\nabla d$, in order to obtain control of
$(\displaystyle\|\nabla\rho\|_{L^\infty_t L^2_x}+
\displaystyle\|\nabla u\|_{L^\infty_t L^2_x}+\displaystyle\|\nabla^2 d\|_{L^\infty_t L^2_x})$, we estimate
$(\displaystyle \|\sqrt{\rho}\dot{u}\|_{L^2_tL^2_x}+\displaystyle \|\nabla d_t\|_{L^2_tL^2_x})$ by
combining an energy estimate of the equation (\ref{clc-2}) in terms of the {\it material derivative} $\dot{u}
\equiv u_t+u\cdot\nabla u$ with second order energy estimates of both  (\ref{clc-2}) and (\ref{clc-3}).
See Lemma 2.4.\\ 
(4) We estimate
$(\displaystyle \|\nabla^2 u\|_{L^\infty_tL^2_x}+\displaystyle \|\nabla^3 d\|_{L^\infty_tL^2_x})$ by
combining thrid order estimate estimates of (\ref{clc-2}) and (\ref{clc-3})  with $H^2$-estimate 
of the Lam\'e equation  and $H^3$-estimate of the harmonic map equation. See Lemma 2.6.
\\ (5) Finally, we obtain the estimate
of $\displaystyle \|\nabla \rho\|_{L^\infty_tL^q_x}$ for $3<q\le 6$ in terms of $\displaystyle
\|u\|_{L^2_tD^{2,q}_x}$. To do it, we  employ the elliptic estimate of the equation satisfied by the {\it effective viscous flux} $G\equiv (2\mu+\lambda){\rm{div}} u-P(\rho)$ and the bound of
$\displaystyle \|\nabla^4 d\|_{L^2_tL^2_x}$ and $\displaystyle\|\nabla u_t\|_{L^2_t L^2_x}$.
See Lemma 2.7.

We would like to point out that during all these arguments, specific forms of the pressume function $P$
play no roles, it is the local Lipschitz regularity of $P$ that is relevant.\\

\noindent{\bf Acknowledgement}. The first two authors are partially supported by NSF grant 1000115.
The work was completed during the third author's visit to the University of Kentucky, which is partially supported by
the second author's NSF grant 0601162. The third author would like to thank the department of Mathematics for its hospitality.

\section {Proof of Theorem \ref{umaintheorem}}
\setcounter{equation}{0}
\setcounter{theorem}{0}

Let $0<T_*<\infty$ be the maximum time for the existence of strong solution $(\rho, u,d)$
to (\ref{clc-1})-(\ref{clc-3}). Namely,  $(\rho, u, d)$ is a strong
solution to (\ref{clc-1})-(\ref{clc-3}) in $\Omega\times (0, T]$ for any $0<T<T_*$, but not
a strong solution in $\Omega\times (0, T_*]$.  Suppose that (\ref{clcblpcondition}) were false,
i.e.
\beq\label{2.1}
M_0:=\int_0^{T_*} \left(\|\mathcal{D}(u)\|_{L^{\infty}}+\|\nabla
d\|^2_{L^{\infty}}\right)\,dt<\infty. \eeq
The goal is to show that under the assumption (\ref{2.1}),
there is a bound $C>0$ depending only on $M_0, \rho_0, u_0, d_0$, and $T_*$ such that
\beq{}\label{uniform_est1}
\sup_{0\le t<T_*}\left[\max_{r=2, q}(\|\rho\|_{W^{1,r}}+\|\rho_t\|_{L^r})
+(\|\sqrt{\rho}u_t\|_{L^2}+\|\nabla u\|_{H^1})+(\|d_t\|_{H^1}+\|\nabla d\|_{H^2})\right]\le C,
\eeq
and
\beq\label{uniform_est2}
\int_{0}^{T_*}\left(\|u_t\|_{D^1}^2+\|u\|_{D^{2,q}}^2+\|d_t\|_{H^2}^2
+\|\nabla d\|_{H^3}^2\right)\,dt\le C.
\eeq
With (\ref{uniform_est1}) and (\ref{uniform_est2}), we can then show without much difficulty
that $T_*$ is not the maximum time, which is the desired contradiction.

Throughout the rest of the paper, we denote by $C$ a generic constant depending only on
$\rho_0$, $u_0$, $d_0$, $T_*$, $M_0$, $\lambda$, $\mu$, $\Omega$, and $P$. We denote by
$$A\lesssim B$$ if there exists a generic constant $C$ such that $A\leq C B$.
For two $3\times 3$ matrices $M=(M_{ij}), N=(N_{ij})$, denote the scalar product  between
$M$ and $N$ by
$$M:N=\sum_{i,j=1}^3 M_{ij} N_{ij}.$$
For $d:\Omega\to S^2$, denote by $\nabla d\otimes \nabla d$ as the $3\times 3$ matrix given by
$$(\nabla d\otimes\nabla d)_{ij}=\langle \nabla_i d, \nabla_j d\rangle, \ 1\le i, j\le 3.$$

The proof is divided into several steps, and we proceed as follows.\\

\noindent{\bf Step 1}. We will first establish $L^\infty$-control of $\rho$. More precisely,  we have

\blm\label{uclemma2.2} 
Let $0<T_*<+\infty$ be the maximum time for a strong
solution $(\rho,u,d)$ to (\ref{clc-1})-(\ref{clc-3}), (\ref{clcinitial}) together with (\ref{clcboundary1})
or (\ref{clcboundary2}) or (\ref{clcboundary3}). If (\ref{first3.1}) and (\ref{2.1}) hold, then 
\beq\label{ucblp2.1} \sup\ls_{0\leq t<T_*}\|\rho\|_{L^{\infty}}\leq C. \eeq\elm

\pf This estimate is a well-known fact of (\ref{clc-1}) and was proved by Huang-Li-Xin \cite{HLX} (Lemma 2.1). 
For the convenience of reader, we sketch it here.
For any $1<r<+\infty$, multiplying (\ref{clc-1}) by $r\rho^{r-1}$ 
and integrating over $\om$, we obtain
\bex
\begin{split}
\frac{d}{dt}\int\rho^{r}\,dx
=&-\int\left(u\cdot\nabla(\rho^{r})+r\rho^{r}\di u\right)\,dx\\
=&-\int\left(\di(u\rho^{r})+(r-1)\rho^{r}\di u\right)\,dx
\leq(r-1)\|\di u\|_{L^{\infty}}\int\rho^{r}\,dx.
\end{split}
\eex
Thus
$$\frac{d}{dt}\|\rho\|_{L^r}\leq \frac{r-1}{r}\|\di u\|_{L^{\infty}}\|\rho\|_{L^r}.$$
This, (\ref{2.1}), Lemma 2.1, together with Gronwall's inequality, imply
$$\sup_{0\le t<T_*}\|\rho(t)\|_{L^r}
\leq \|\rho_0\|_{L^r}\exp\left(\int_0^{T_*}\|\di u\|_{L^{\infty}}dt\right)\le C,$$
which, after sending $r$ to $\infty$, implies (\ref{ucblp2.1}).
This completes the proof.
\endpf\\

\noindent{\bf Step 2}. We next establish the global energy inequality for strong solutions, namely,

\blm\label{uclemma2.1}{\it Let $0<T_*<+\infty$ be the maximum time for a strong
solution $(\rho,u,d)$ to (\ref{clc-1})-(\ref{clc-3}), (\ref{clcinitial}) together with (\ref{clcboundary1})
or (\ref{clcboundary2}) or (\ref{clcboundary3}). If (\ref{first3.1}) and (\ref{2.1}) hold, then for any $0\le t<T_*$,
the following inequality holds:
\beq\label{clcblp2.1}
\begin{split}
&\int_\Omega\left(\rho|u|^2+|\nabla d|^2\right)(t)\,dx
+\int_0^t\int_{\om}\Big(|\nabla u|^2+\left|\de d+|\nabla d|^2d\right|^2\Big)\,dx\,ds\\
&\le C\Big[\int_{\om}\left(\rho_0|u_0|^2+|\nabla d_0|^2\right)\,dx+1\Big].
\end{split}
\eeq 
Furthermore, we have
\beq\label{clcblp2.2}
\int_0^{T_*}\int_{\om}|\nabla^2 d|^2\,dx\,dt\leq C. \eeq}
 \elm

\pf Without loss of generality, we may assume that $P(0)=0$. Since $P$ is locally Lipschitz by
(\ref{regularity_p}), it follows that $P'$ is locally bounded  on $[0,+\infty)$.
Since $\rho$ is bounded in $\Omega\times [0, T_*)$ by (\ref{ucblp2.1}), we then have that, on
$\Omega\times [0,T_*)$, 
\begin{eqnarray}\label{upper_bound_p} 
|P(\rho)|&\le& \left\|P'(\rho)\right\|_{L^\infty}\rho\le C\rho\ (\ \le C\ )\\
|\nabla(P(\rho))|&\le& \left\|P'(\rho)\right\|_{L^\infty}|\nabla\rho|\le C|\nabla\rho|.\nonumber
\end{eqnarray}
Since (\ref{clcboundary1}) and (\ref{clcboundary2}) are easier to handle\footnote{in fact, with respect to the boundary condition (\ref{clcboundary1}) and (\ref{clcboundary2}), by integration by parts one has
$$-\int (\mu\Delta u+(\mu+\lambda)\nabla({\rm{div}} u))\cdot u\,dx
=\int (\mu|\nabla u|^2+(\mu+\lambda)|{\rm{div}} u|^2)\,dx\ge \mu \int |\nabla u|^2\,dx.$$},
we outline the proof for the boundary condition (\ref{clcboundary3}).
 Multiplying (\ref{clc-2}) by $u$ and integrating over $\om$, we have
\begin{eqnarray}\label{clcblp2.3}
&&\frac{1}{2}\int \rho(\partial_t|u|^2+u\cdot\nabla |u|^2)\,dx+\int (\mu|\nabla \times u|^2+(2\mu+\lambda)|\mbox{div}u|^2)\,dx\nonumber\\
&=&\int P(\rho)\mbox{div}u\, dx-\int u\cdot\nabla d\cdot \de d\,dx.
\end{eqnarray}
Here we have used the fact that $\de u=\nabla \mbox{div} u-\nabla\times\curl u$, and the Navier-slip boundary condition (\ref{clcboundary3})
to obtain
\bex
\int \mathcal L u\cdot u\,dx
&=&\int[(2\mu+\lambda)\nabla(\mbox{div}u)\cdot u -\mu\nabla\times\curl u\cdot u]\,dx\\
&=&-\int  (2\mu+\lambda)|\mbox{div}u|^2+\mu |\curl u|^2)\,dx.\\
\eex
Hereafter we repeatedly use the following identity:
$$\int \langle\nabla\times u, \curl u\rangle\,dx
=\int\langle u, \nabla\times(\curl u)\rangle\,dx, \ \forall u \ {\rm{with}}\  \curl u\times\nu=0 \ {\rm{on}}\ \pa\om.$$
By the formula of transportion, we have
$$\int \rho(\partial_t|u|^2+u\cdot\nabla |u|^2)\,dx=\frac{d}{dt}\int\rho|u|^2\,dx.$$
By (\ref{upper_bound_p}) and Cauchy's inequality, we have
\begin{eqnarray*}
|\int P(\rho)\mbox{div}u\,dx|&\leq &\int |P(\rho)||{\rm{div}} u|\,dx
\lesssim \int \rho|\nabla u|\,dx\\
&\le&\epsilon\int |\nabla u|^2\,dx+C(\epsilon)\int \rho^2
\le \epsilon\int |\nabla u|^2\,dx+C(\epsilon),
\end{eqnarray*}
where we have used (\ref{ucblp2.1}) and the conservation of mass to assure 
$$\int \rho^2 \le \|\rho\|_{L^1}\|\rho\|_{L^\infty}\le C.$$
Putting these inequalities into (\ref{clcblp2.3}), we obtain
\beq\label{clcblp2.4}
\frac{d}{dt}\int \rho|u|^2\,dx
+\int(\mu|\nabla\times u|^2+(2\mu+\lambda)|\mbox{div}u|^2)\,dx\le -\int u\cdot\nabla d\cdot\de d\,dx
+\epsilon\int |\nabla u|^2\,dx+C(\epsilon).
\eeq
Since $\Omega$ is assumed to be simply connected for the boundary condition (\ref{clcboundary3}),
we have the following estimate (see \cite{Wahl} for its proof):
\beq\label{div-curl0}
\|\nabla u\|_{L^2}\les \|\nabla\times u\|_{L^2}+\|\di u\|_{L^2} \ \forall u\in H^1(\Omega) \ {\rm{with}}\ u\cdot\nu=0 \ {\rm{on}}\ \pa\om.
\eeq
This, combined with (\ref{viscosity}), implies that
\beq\label{div-curl} \int(\mu|\nabla\times u|^2+(2\mu+\lambda)|\mbox{div}u|^2)dx\ge
\frac{\mu}{3}\int (|\nabla \times u|^2+|{\rm{div}} u|^2)\,dx\ge \frac{1}{C}\int |\nabla u|^2\,dx.
\eeq
Thus, by choosing $\epsilon=\frac{1}{2 C}$, (\ref{clcblp2.4}) implies
\beq\label{clcblp2.40}
\frac{d}{dt}\int \rho|u|^2\,dx
+\frac{1}{2C}\int |\nabla u|^2\,dx\le -\int u\cdot\nabla d\cdot\de d\,dx
+C.
\eeq
Now, multiplying (\ref{clc-3}) by $(\de d+|\nabla d|^2d)$, integrating over $\om$ and using $\frac{\partial d}{\partial\nu}=0$ on $\partial\Omega$, we have
\beq\label{clcblp2.5}
\frac{1}{2}\frac{d}{dt}\int |\nabla d|^2\,dx+\int\left|\de d+|\nabla d|^2d\right|^2dx=\int u\cdot\nabla d\cdot\de d\,dx,
\eeq
where we have used the fact that $|d|=1$ in $\Omega$ and hence
$$\int \langle d_t+u\cdot\nabla d, |\nabla d|^2d \rangle\,dx=0.$$
Adding (\ref{clcblp2.40}) and (\ref{clcblp2.5}) together yields (\ref{clcblp2.1}).

To see (\ref{clcblp2.2}), observe that (\ref{2.1}) implies  $\int_0^{T_*}\|\nabla d\|_{L^{\infty}}^2\,dt\le M_0$ so that 
\bex
\int_0^{T_*}\int_\om |\nabla d|^4\,dx\,dt&\leq& M_0
\cdot\left(\sup\ls_{0\leq t<T_*}\int |\nabla d|^2\,dx\right)\\
&\leq& CM_0\left[1+\int\left(\rho_0|u_0|^2+|\nabla d_0|^2\right)\,dx\right]
\eex
where we have used (\ref{clcblp2.1}) in the last step. This and (\ref{clcblp2.1}) then imply
\bex 
\int_0^{T_*}\int_{\om}|\de d|^2\,dx\,dt&=&\int_0^{T_*}\int_{\om}\left|\de d+|\nabla d|^2d\right|^2\,dx\,dt
+\int_0^{T_*}\int_{\om}|\nabla d|^4\,dx\,dt\\
&\leq& C M_0\Big[1+\int\left(\rho_0|u_0|^2+|\nabla d_0|^2\right)\,dx\Big].
\eex
Since $\frac{\partial d}{\partial\nu}=0$ on $\partial\Omega$, the standard $L^2$-estimate yields
$$\int |\nabla^2 d|^2\,dx\le C\int (|\de d|^2+|\nabla d|^2)\,dx.$$
Thus (\ref{clcblp2.2}) follows easily, and the proof is complete.
\endpf

\vspace{5mm}

\noindent{\bf Step 3}. We will establish $L^\infty_t L^r_x$-control of $\nabla d$ for any $2\le r<+\infty$,
a key ingredient for the higher order estimates of $u,\nabla d$.  More precisely, we have

\blm\label{uclemma2.3}{\it Let $0<T_*<+\infty$ be the maximum time for a strong
solution $(\rho,u,d)$ to (\ref{clc-1})-(\ref{clc-3}), (\ref{clcinitial}) together with (\ref{clcboundary1})
or (\ref{clcboundary2}) or (\ref{clcboundary3}). If (\ref{first3.1}) and (\ref{2.1}) hold, then 
for any $2\le r<+\infty$, there exists a $C>0$ depending on $M_0, u_0,d_0, \Omega, n,$ and $r$
such that 
\beq\label{ucblp2.2} \sup\ls_{0\leq t<T_*}\|\nabla
d\|_{L^{r}}^r+\int_0^{T_*}\int_{\om}|\nabla d|^{r-2}|\nabla^2
d|^2\,dx\,dt\leq C. \eeq} \elm

\pf Here we only consider the Navier-slip boundary condition (\ref{clcboundary3}), since the argument to deal the first two boundary conditions (\ref{clcboundary1}) and (\ref{clcboundary2}) is similar and easier.
Differentiating the equation (\ref{clc-3}) with rrespect to $x$, we have
\beq\label{ucblp2.3}
\nabla d_t-\nabla\de d=-\nabla(u\cdot\nabla d)+\nabla(|\nabla d|^2d).
\eeq
Multiplying (\ref{ucblp2.3}) by $r|\nabla d|^{r-2}\nabla d$ and integrating over $\om$, we obtain
\beq\label{ucblp2.4}
\begin{split}
&\frac{d}{dt}\int|\nabla d|^{r}\,dx+r\int\left(|\nabla d|^{r-2}|\nabla^2 d|^2+(r-2)|\nabla d|^{r-2}|\nabla |\nabla d||^2\right)\,dx\\
=&r\int\nabla(|\nabla d|^2d)|\nabla d|^{r-2}\nabla d\,dx
-r\int\nabla(u\cdot\nabla d)|\nabla d|^{r-2}\nabla d\,dx\\
&+\frac{r}2\int_{\pa\om}|\nabla d|^{r-2}\langle \nabla(|\nabla d|^2),\nu\rangle\,d\sigma
=\sum\ls_{i=1}^3I_i.
\end{split}
\eeq
We can estimate $I_i$ ($i=1,2,3$) separately as follows.
For $I_1$, since $$\nabla(|\nabla d|^2 d)=|\nabla d|^2\nabla d+\nabla(|\nabla d|^2) d
\mbox{ and }d\cdot\nabla d=0,$$
we have  
\bex\begin{split}
I_1=&r\int|\nabla d|^{r+2}\,dx
\les\|\nabla d\|_{L^{\infty}}^2\int|\nabla d|^{r}\,dx.
\end{split}
\eex
For $I_2$, we have
\bex\begin{split}
I_2=&-r\int|\nabla d|^{r-2}\nabla_i u^j\langle\nabla_j d,\nabla_i d\rangle\,dx
-\int u\cdot\nabla(|\nabla d|^{r})\,dx\\
=&-r\int|\nabla d|^{r-2}\mathcal{D} (u):\nabla d\otimes\nabla d\,dx+
\int(\di u)|\nabla d|^{r}\,dx
\les\|\mathcal{D}(u)\|_{L^{\infty}}\int|\nabla d|^{r}\,dx.
\end{split}
\eex
The estimate of the boundary integral $I_3$ is more delicate.  
Let $\mathbb {I}_{\partial\om}$ denote the second fundamental form of $\pa\om$: for any $x\in\pa\om$, 
$$\mathbb I_{\pa\om}(x)(U,V)=-\nabla\nu(x)(U,V), \ \forall U, V\in T_x(\pa\om).$$
Let $\nabla_T$ denote the tangential derivative on $\partial\om$. Since 
$\frac{\pa d}{\pa\nu}=\langle\nabla d,\nu\rangle=0 \ \mbox{on}\ \pa\om, $
we have $\nabla_T(\frac{\pa d}{\pa\nu})=0  \ \mbox{on}\ \pa\om$. Hence we have, on $\partial\Omega$, 
$$\frac12\langle\nabla(|\nabla d|^2),\nu\rangle
=\nabla d\cdot \nabla\langle \nabla d,\nu\rangle-\nabla \nu(\nabla d,\nabla d)
=\nabla_T d\cdot \nabla_T(\frac{\pa d}{\pa \nu})-\nabla\nu(\nabla_T d,\nabla_T d)
=\mathbb I_{\pa\om}(\nabla_T d,\nabla_Td).$$
Therefore we have
\bex\begin{split}
I_3=&r\int_{\pa\om}|\nabla d|^{r-2}\mathbb I_{\pa\om}(\nabla_T d,\nabla_Td)\,d\sigma
\les \int_{\pa\om}|\nabla d|^r\,d\sigma.
\end{split}
\eex
Applying the trace formula $W^{1,1}(\om)\subset L^1(\pa\om)$ and H\"older's inequality, we obtain
\bex
\begin{split}
I_3 \les& \||\nabla d|^r\|_{W^{1,1}}
\les\int|\nabla d|^r\,dx+\int|\nabla d|^{r-1}|\nabla^2 d|\,dx\\
\le& C\int|\nabla d|^r\,dx+\frac{r}{4}\int|\nabla d|^{r-2}|\nabla^2 d|^2\,dx.
\end{split}
\eex
Putting all these estimates into (\ref{ucblp2.4}), we obtain
\bex
\begin{split}
\frac{d}{dt}\int|\nabla d|^{r}\,dx+\frac{r}{2}\int|\nabla d|^{r-2}|\nabla^2 d|^2\,dx
\les\left(\|\nabla d\|_{L^{\infty}}^2+\|\mathcal{D}(u)\|_{L^{\infty}}+1\right)\int|\nabla d|^{r}\,dx.
\end{split}
\eex
By Gronwall's inequality and (\ref{2.1}),  we obtain that for 
any $0\le t<T_*$,
\bex
\begin{split}
&\int|\nabla d(t)|^{r}\,dx+\int_0^t\int_{\om}|\nabla d|^{r-2}|\nabla^2 d|^2\,dxds\\
\les&\int|\nabla d_0|^{r}\,dx\cdot\exp\left(\int_0^{T_*}
(1+\|\nabla d\|_{L^{\infty}}^2+\|\mathcal{D}(u)\|_{L^{\infty}})\,dt\right)
\le \ C.
\end{split}
\eex This completes the proof.
\endpf\\

\noindent{\bf Step 4}. Estimates of $(\nabla u,\nabla \rho,\nabla^2 d)$ in $L^\infty_tL^2_x(\om\times [0,T_*])$.  First, for any function $f$ on $\Omega\times (0,T_*)$,
let $$\dot{f}=f_t+u\cdot\nabla f$$
denote the material derivative of $f$.
Then we have

\blm\label{uclemma2.4}{\it Let $0<T_*<+\infty$ be the maximum time for a strong
solution $(\rho,u,d)$ to (\ref{clc-1})-(\ref{clc-3}), (\ref{clcinitial}) together with (\ref{clcboundary1})
or (\ref{clcboundary2}) or (\ref{clcboundary3}). If (\ref{first3.1}) and (\ref{2.1}) hold, then 
\beq\label{ucblp2.5}
\begin{split}
\sup\ls_{0\leq t<T_*}\left(\|\nabla u\|_{L^{2}}^2+\|\nabla
\rho\|_{L^{2}}^2+\|\nabla^2 d\|_{L^{2}}^2\right)
+\int_0^{T_*}\int_{\om}\left(\rho|\dot{u}|^{2}+|\nabla
d_t|^2\right)\,dxdt\leq C.
\end{split}
\eeq} \elm

\pf To make the presentation shorter, here we only consider the difficult case: the Navier-slip boundary condition (\ref{clcboundary3}). To obtain
the estimates of $u$, we adapt some arguments by \cite{HLX} (Lemma 2.2). 
Multiplying (\ref{clc-2}) by $\dot{u}$ and integrating over $\om$, we obtain,
\beq\label{ucblp2.6}
\begin{split}
&\int\rho|\dot{u}|^2\,dx
-\int\langle\mathcal Lu, u_t\rangle\,dx\\
=&\int\langle u\cdot\nabla u,\mathcal Lu\rangle\,dx-\int u\cdot\nabla u\cdot\nabla (P(\rho))\,dx-\int u_t\cdot\nabla (P(\rho))\,dx\\
-&\int u\cdot\nabla u\cdot\langle\de d,\nabla d\rangle\,dx
-\int u_t\cdot\langle\de d,\nabla d\rangle\,dx.
\end{split}
\eeq
Similar to the proof of Lemma 2.1, we have
\bex\begin{split}
-\int\langle\mathcal Lu, u_t\rangle\,dx
=&\int[\mu\langle\nabla\times\curl u, u_t\rangle-
(2\mu+\lambda)\langle\nabla({\rm{div}}u), u_t\rangle]\,dx\\
=&\int[\mu\langle \nabla\times u, \nabla\times u_t\rangle
+(2\mu+\lambda)({\rm{div}}u)({\rm{div}}u_t)]\,dx\\
=&\frac12\frac{d}{dt}\int[\mu|\nabla\times u|^2+(2\mu+\lambda)({\rm{div}}u)^2]\,dx,\\
\end{split}
\eex
where we have used the fact that $u_t\cdot\nu=\curl u\times\nu=0$ on $\pa\om$ during
the integration by parts.

The terms on the right hand side of (\ref{ucblp2.6}) can be estimated as follows.
\beq\label{ucblp2.7}
\begin{split}
\int\langle u\cdot\nabla u, \mathcal Lu\rangle\,dx
=&-\mu\int \langle u\cdot\nabla u, \nabla\times\curl  u\rangle\,dx
+(2\mu+\lambda)\int\langle u\cdot\nabla u,\nabla(\di u)\rangle\,dx.
\end{split}
\eeq
For the first term in the right hand side of (\ref{ucblp2.7}), by using $\curl u\times\nu=0$
on $\pa\om$ and the formula $$u\times\curl u=\frac12\nabla(|u|^2)-u\cdot\nabla u,$$
we have
\bex
\begin{split}
&-\mu\int \langle u\cdot\nabla u, \nabla\times\curl u\rangle\,dx=
-\mu\int \langle\curl u,\nabla\times (u\cdot\nabla u)\rangle\,dx
=\mu\int \curl u\cdot\nabla\times (u\times\curl u) \,dx\\
=&\mu\int\langle\curl u,
(\curl u\cdot\nabla) u-(u\cdot\nabla) \curl u+\di(\curl u)u-(\di u)\curl u\rangle \,dx\\
=&\mu\int\left(\mathcal D(u): \curl u\otimes \curl u-\frac{1}{2}\di u(\curl u)^2\right) \,dx
\les\|\mathcal{D}(u)\|_{L^{\infty}}\|\nabla u\|_{L^2}^2,
\end{split}
\eex
where we have also used the formulas
$$\nabla\times(a\times b)=(b\cdot\nabla) a-(a\cdot\nabla) b+({\rm{div}}b)a-
({\rm{div}}a)b,\ {\rm{ and }}  \ {\rm{div}}(\curl u)=0.$$

To estimate the second term in the right hand side  of (\ref{ucblp2.7}),
denote by $u^\tau=u-(u\cdot\nu)\nu$ the tangential component of $u$ on $\pa\om$.
Note that (\ref{clcboundary3}) implies $u=u^\tau$ on $\pa\om$. Hence 
we have
\bex\begin{split}
\int\langle u\cdot\nabla u, \nabla(\di u)\rangle\,dx
=&\int_{\pa\om}\langle (u\cdot\nabla) u, \nu\rangle\di u\,d\sigma
-\int ((\nabla u):(\nabla u)^t\di u+\frac12u\cdot\nabla((\di u)^2))\,dx\\
=&\int_{\pa\om}u^\tau\cdot\nabla_T (u\cdot\nu)\di u\,d\sigma
-\int_{\pa\om}\nabla \nu(u^\tau, u^\tau)\di u\,d\sigma\\
-&\int[(\nabla u):(\nabla u)^t\di u-\frac12(\di u)^3]\,dx\\
=&\int_{\pa\om}\mathbb I_{\pa\om}(u^\tau,u^\tau)\di u\,d\sigma
-\int[(\nabla u):(\nabla u)^t\di u-\frac12(\di u)^3]\,dx\\
\les&\|u\|_{L^{4}(\pa\om)}^2\|\di u\|_{L^2(\pa\om)}+\|\mathcal{D}(u)\|_{L^{\infty}}\|\nabla u\|_{L^2}^2.
\end{split}
\eex
By the trace formula $H^1(\om)\subset L^r(\partial\om)$ for $r=2,4$, the Poincar\'e inequality (see \cite{yoshida-giga}):
$$\|u\|_{L^2}\les\|\nabla u\|_{L^2}, \ \forall u\in H^1 \ \mbox{with}\ u\cdot\nu=0 
\ \mbox{on}\ \pa\om,$$
and H\"older's inequality, we have
\bex\begin{split}
\|u\|_{L^{4}(\pa\om)}^2\|\di u\|_{L^2(\pa\om)}
\les &\|u\|_{H^1}^2\|\nabla u\|_{H^1}
\les  \|\nabla u\|_{L^2}^2(\|\nabla u\|_{L^2}+\|\nabla^2 u\|_{L^2})\\
\leq &C(\epsilon)(1+\|\nabla u\|_{L^2}^4)+\epsilon\|\nabla^2 u\|_{L^2}^2
\end{split}
\eex
for small $\epsilon>0$ to be determined later. Thus we obtain 
\beq\label{ucblp2.8}
\begin{split}
\int\langle u\cdot\nabla u, \mathcal Lu\rangle\,dx
\leq\epsilon\|\nabla^2 u\|_{L^2}^2 
+C(\epsilon)(1+\|\nabla u\|_{L^{2}}^4)+C\|\mathcal{D}(u)\|_{L^{\infty}}\|\nabla u\|_{L^2}^2.
\end{split}
\eeq
The remaining terms in the right hand side of (\ref{ucblp2.6}) can be estimated as follows.
\beq\label{ucblp2.9}
\begin{split}
-&\int u\cdot\nabla u\cdot\nabla (P(\rho))\,dx\\
=&-\int_{\pa\om}P(\rho) \langle (u\cdot\nabla) u, \nu\rangle\,d\sigma
+\int (P(\rho) u\cdot\nabla (\di u)+P(\rho)(\nabla u):(\nabla u)^t)\,dx\\
=&-\int_{\pa\om}[P(\rho)(u\cdot\nabla) (u\cdot\nu)-P(\rho)\mathbb I_{\pa\om}(u^\tau,u^\tau)]\,d\sigma
+\int P(\rho) (\nabla u):(\nabla u)^t\,dx\\
&+\int_{\pa\om}P(\rho)(u\cdot\nu)\di u\,d\sigma
-\int[\nabla (P(\rho))\cdot u\di u+P(\rho)(\di u)^2)]\,dx\\
=&\int_{\pa\om}P(\rho)\mathbb I_{\pa\om}(u^\tau,u^\tau)\,d\sigma
+\int [P(\rho) ((\nabla u):(\nabla u)^t-({\rm{div}}u)^2)]\,dx\\
&-\int\nabla (P(\rho))\cdot u\di u\,dx\\
\les&\| u\|_{L^{2}(\pa\om)}^2+\|\nabla u\|_{L^{2}}^2+\int_{\om}|\nabla \rho||u||\di u|\,dx\\
\les&\|\nabla u\|_{L^{2}}^2+\|\mathcal{D} (u)\|_{L^{3}}\| u\|_{L^{6}}\|\nabla \rho\|_{L^{2}}\\
\les&\|\nabla u\|_{L^{2}}^2+\|\mathcal{D}(u)\|_{L^{\infty}}^{\frac{1}{3}}\|\nabla u\|_{L^{2}}^{\frac{5}{3}}\|\nabla \rho\|_{L^{2}}\\
\les&1+\left(\|\mathcal{D} (u)\|_{L^{\infty}}+1\right)\|\nabla u\|_{L^{2}}^2+\|\nabla u\|_{L^{2}}^2\|\nabla \rho\|_{L^{2}}^2,
\end{split}
\eeq
where we have used (\ref{upper_bound_p}) and the Sobolev and Poincar\'e inequalities (see \cite{yoshida-giga}):
$$\|u\|_{L^6}\les (\|u\|_{L^2}+\|\nabla u\|_{L^2})\les\|\nabla u\|_{L^2},
\ \forall u \in H^1(\Omega) \ {\rm{with}}\ u\cdot\nu =0 \ {\rm{on}}\ \pa\om.$$
Since (\ref{ucblp2.1}) and (\ref{clc-1}) also imply
\beq\label{upper_bound_p_t} |(P(\rho))_t|
\le \|P'(\rho)\|_{L^\infty} |\rho_t|\le\|P'(\rho)\|_{L^\infty} (\rho|\nabla u|+|\nabla\rho||u|)\le C(|\nabla u|+|u||\nabla\rho|),
\eeq
we have
\beq\label{ucblp2.10}
\begin{split}
-&\int u_t\cdot\nabla P(\rho)\,dx=\frac{d}{dt}\int P(\rho)\di u\,dx-\int (P(\rho))_t\di u\,dx\\
\le &\frac{d}{dt}\int P(\rho)\di u\,dx+C\int (|u||\nabla \rho| |{\rm{div}} u|+|\nabla u|^2)\,dx\\
\le &\frac{d}{dt}\int P(\rho)\di u\,dx
+C\left[\|\nabla u\|_{L^{2}}^2+\|\mathcal D(u)\|_{L^{3}}\|u\|_{L^6}\|\nabla \rho\|_{L^{2}}\right]\\
\le & \frac{d}{dt}\int P(\rho)\di u\,dx
+C\left[\|\nabla u\|_{L^{2}}^2+\|\mathcal D(u)\|_{L^\infty}^\frac13\|\nabla u\|_{L^2}^\frac53\|\nabla \rho\|_{L^{2}}\right]\\
\le &  \frac{d}{dt}\int P(\rho)\di u\,dx
+C\left[1+(1+\|\mathcal D(u)\|_{L^\infty})\|\nabla u\|_{L^2}^2+\|\nabla u\|_{L^2}^2\|\nabla \rho\|_{L^{2}}^2\right],\\
\end{split}
\eeq
where we have used the Poincar\'e inequality and H\"older's inequality.
\beq\label{ucblp2.11}
\begin{split}
-\int u\cdot\nabla u\cdot\langle\de d,\nabla d\rangle\,dx
\leq&\|u\|_{L^{6}}\|\nabla u\|_{L^{6}}\|\de d\|_{L^{2}}\|\nabla d\|_{L^{6}}\\
\les&\|\nabla u\|_{L^{2}}(\|\nabla u\|_{L^2}+\|\nabla^2 u\|_{L^{2}})\|\de d\|_{L^{2}} \ \ \ ( {\rm{by}}\ (\ref{ucblp2.2}) 
\ {\rm{with}}\ r=6)\\
\le& \epsilon\|\nabla^2 u\|_{L^{2}}^2+C(\epsilon)\|\nabla u\|_{L^{2}}^2\|\de d\|_{L^{2}}^2+
\|\nabla u\|_{L^2}^4 +\|\de d\|_{L^2}^2.
\end{split}
\eeq
To estimate the last term in the right hand side of (\ref{ucblp2.6}), 
denote $M(d)=\nabla d\otimes\nabla d-\frac12|\nabla d|^2\mathbb I_3$. 
Then we have $\langle\de d,\nabla d\rangle ={\rm{div}}(M(d))$ and 
\beq\label{ucblp2.12}
\begin{split}
-\int u_t\cdot\langle\de d,\nabla d\rangle\,dx=&-\int_{\pa\om}u_t\cdot (\langle M(d),\nu\rangle)\,d\sigma+\int M(d):\nabla u_t\,dx\\
=&\int M(d):\nabla u_t\,dx\ \ ({\rm{since}}\ u_t\cdot\nu=\frac{\pa d}{\pa\nu}=0\ {\rm{on}}\ \pa\om)\\
=&\frac{d}{dt}\int M(d):\nabla u\,dx-\int (M(d))_t:\nabla u\,dx\\
\le&\frac{d}{dt}\int M(d):\nabla u\,dx+C\int |\nabla d_t||\nabla d||\nabla u|\,dx\\
\le&\frac{d}{dt}\int M(d):\nabla u\,dx+C(\epsilon)
\|\nabla d\|_{L^{\infty}}^2\|\nabla u\|_{L^{2}}^2+\epsilon\|\nabla d_t\|_{L^{2}}^2.
\end{split}
\eeq
Putting (\ref{ucblp2.8})-(\ref{ucblp2.12}) into (\ref{ucblp2.6}), we obtain
\beq\label{ucblp2.13}
\begin{split}
&\frac{1}{2}\frac{d}{dt}\int
\left(\mu|\nabla \times u|^{2}+(2\mu+\lambda)|\di u|^{2}\right)\,dx+\int
\rho|\dot{u}|^2\,dx\\
\le &\frac{d}{dt}\int (M(d):\nabla u+P(\rho)\di u)\,dx+\epsilon(\|\nabla d_t\|_{L^{2}}^2+\|\nabla^2u\|_{L^{2}}^2)\\
&+C\left(\|\mathcal{D} (u)\|_{L^{\infty}}+\|\nabla u\|_{L^{2}}^2+\|\nabla d\|_{L^{\infty}}^2+1\right)\|\nabla u\|_{L^{2}}^2\\
&+C\|\nabla u\|_{L^{2}}^2\|\nabla \rho\|_{L^{2}}^2
+C(1+\|\nabla u\|_{L^{2}}^2)\|\de d\|_{L^{2}}^2+C(\epsilon).
\end{split}
\eeq

Now we want to estimate $\|\nabla\rho\|_{L^2}^2$. 
Differentiating the equation (\ref{clc-1}) with respect to $x$, 
multiplying the resulting equation by $2\nabla\rho$ and integrating over $\om$, we obtain
\beq\label{ucblp2.14}
\begin{split}
\frac{d}{dt}\|\nabla \rho\|_{L^{2}}^2=&-\int (u\cdot\nabla(|\nabla\rho|^2)+2\mathcal D(u):\nabla\rho\otimes\nabla\rho
+2|\nabla \rho|^2(\di u)+2\rho\nabla(\di u)\cdot\nabla\rho)\,dx\\
=&-\int |\nabla \rho|^2\di u\,dx-2\int \mathcal{D}(u):\nabla\rho\otimes\nabla\rho\,dx-2\int
\rho\nabla\di u\cdot\nabla\rho\,dx\\
\les&\left(\|\mathcal{D}(u)\|_{L^{\infty}}+1\right)\|\nabla \rho\|_{L^2}^2
+\epsilon\|\nabla^2 u\|_{L^2}^2.
\end{split}
\eeq

Next we want to estimate $\|\nabla d_t\|_{L^2}^2$. To do this, we 
multiply (\ref{ucblp2.3}) by $\nabla d_t$ and integrate over $\om$ and use
$\frac{\partial d_t}{\partial\nu}=0$ on $\pa\om$ to obtain
\beq\label{ucblp2.15}
\begin{split}
&\frac{d}{dt}\int |\de d|^2dx+\int |\nabla d_t|^2dx
=\int \left(\nabla(|\nabla d|^2d)-\nabla(u\cdot\nabla d)\right)\nabla d_tdx\\
\leq&C(\epsilon)\int \left(|\nabla u|^2|\nabla d|^2+|u|^2|\nabla^2 d|^2+|\nabla d|^6+|\nabla d|^2|\nabla^2 d|^2\right)dx+\epsilon\|\nabla d_t\|_{L^2}^2\\
\le&\epsilon\|\nabla d_t\|_{L^2}^2+C+C\|\nabla d\|_{L^{\infty}}^2\left(\|\nabla^2 d\|_{L^2}^2+\|\nabla u\|_{L^2}^2\right)+C\int | u|^2|\nabla^2 d|^2dx,
\end{split}
\eeq
where we have used (\ref{ucblp2.2}) (with $r=6$) in the last step.
For the last term in the right hand side of (\ref{ucblp2.15}), by using Nirenberg's interpolation inequality
and (\ref{ucblp2.2}) we have
\begin{equation}\begin{split}
\int |u|^2|\nabla^2 d|^2dx\le&\|u\|_{L^6}^{2}\|\nabla^2 d\|_{L^{3}}^2
\les\|\nabla u\|_{L^2}^2\|\nabla d\|_{L^6}\|\nabla^3 d\|_{L^2}+\|\nabla u\|_{L^2}^2\\
\les&\|\nabla u\|_{L^2}^2\|\nabla^3 d\|_{L^{2}}+\|\nabla u\|_{L^2}^2
\le C(\epsilon)\|\nabla u\|_{L^2}^4+\epsilon\|\nabla^3 d\|_{L^{2}}^{2}+C. \label{cross-est1}
\end{split}\end{equation}
Applying the standard $H^3$-estimate to the Neumann boundary value problem of
the equation (\ref{ucblp2.3}), and using (\ref{ucblp2.2}), we have
\begin{equation}\begin{split}
\|\nabla^3 d\|_{L^2}^{2}\les& \|\nabla\de d\|_{L^2}^2+\|\nabla d\|_{H^1}^2
\les\|\nabla d_t\|_{L^2}^{2}+\|\nabla(u\cdot\nabla d)\|_{L^2}^{2}+\|\nabla( |\nabla d|^2d)\|_{L^2}^{2}
+\|\nabla d\|_{H^1}^2\\
\les&\|\nabla d_t\|_{L^2}^{2}+\|\nabla d\|_{L^{\infty}}^{2}\left(\|\nabla u\|_{L^2}^{2}+\|\nabla^2 d\|_{L^2}^{2}\right)\\
&+\|\nabla d\|_{L^6}^{6}+\|\nabla d\|_{H^1}^2+\int |u|^2|\nabla^2 d|^2dx\\
\les&\|\nabla d_t\|_{L^2}^{2}+\|\nabla d\|_{L^{\infty}}^{2}\left(\|\nabla u\|_{L^2}^{2}+\|\nabla^2 d\|_{L^2}^{2}\right)+\int |u|^2|\nabla^2 d|^2\,dx+
\|\nabla^2 d\|_{L^2}^2+1. \label{cross-est2}
\end{split}\end{equation}
Substituting (\ref{cross-est1}) into (\ref{cross-est2}) and choosing $\epsilon$ sufficiently small, we have
\beq\begin{split}
\|\nabla^3 d\|_{L^2}^{2}
\les 1+\|\nabla d_t\|_{L^2}^{2}+\|\nabla d\|_{L^{\infty}}^{2}\left(\|\nabla u\|_{L^2}^{2}+\|\nabla^2 d\|_{L^2}^{2}\right)+\|\nabla u\|_{L^2}^4
+\|\nabla^2 d\|_{L^2}^2. \label{H3-est}
\end{split}
\eeq
Substituting (\ref{H3-est}) into (\ref{cross-est1}), we obtain
\begin{equation}\label{cross-est3}
\begin{split}
\int |u|^2|\nabla^2 d|^2
\leq C+\epsilon\|\nabla d_t\|_{L^2}^{2}+C\|\nabla d\|_{L^{\infty}}^{2}\left(\|\nabla u\|_{L^2}^{2}+\|\nabla^2 d\|_{L^2}^{2}\right)+C\|\nabla u\|_{L^2}^4.
\end{split}
\end{equation}
Putting (\ref{cross-est3}) into (\ref{ucblp2.15}) and choosing $\epsilon$ sufficiently small, we obtain
\beq\label{ucblp2.16}
\begin{split}
\frac{d}{dt}\int |\de d|^2\,dx+\int |\nabla d_t|^2\,dx
\les&1+\|\nabla u\|_{L^2}^4+\|\nabla d\|_{L^{\infty}}^2\left(\|\nabla^2 d\|_{L^2}^2+\|\nabla u\|_{L^2}^2\right)+\|\nabla^2 d\|_{L^2}^2\\
\les&1+\|\nabla u\|_{L^2}^4+\|\nabla d\|_{L^{\infty}}^2\left(\|\nabla d\|_{L^2}^2+\|\de d\|_{L^2}^2+\|\nabla u\|_{L^2}^2\right)\\
&+\|\Delta d\|_{L^2}^2.
\end{split}
\eeq
Putting (\ref{ucblp2.13}), (\ref{ucblp2.14}) and (\ref{ucblp2.16}) together, we obtain
\bex
\begin{split}
&\frac{d}{dt}\int \left(\mu|\nabla\times u|^{2}+(2\mu+\lambda)|\di u|^{2}+|\nabla\rho|^2+|\de d|^2\right)\,dx
+\int \left(2\rho|\dot{u}|^2+|\nabla d_t|^2\right)\,dx\\
\leq&2\frac{d}{dt}\int(M(d):\nabla u+P(\rho)\di u)\,dx+
\epsilon\|\nabla d_t\|_{L^{2}}^2+\epsilon\|\nabla^2 u\|_{L^{2}}^2\\
+&C(1+\|\nabla u\|_{L^{2}}^2+\|\mathcal D(u)\|_{L^\infty})\|\nabla \rho\|_{L^{2}}^2
+C\left(1+\|\mathcal{D} (u)\|_{L^{\infty}}+\|\nabla u\|_{L^{2}}^2+\|\nabla d\|_{L^{\infty}}^2\right)\|\nabla u\|_{L^{2}}^2\\
+&C(1+\|\nabla u\|_{L^{2}}^2+\|\nabla d\|_{L^{\infty}}^2)\|\de d\|_{L^{2}}^2+C\|\nabla d\|_{L^\infty}^2+C.
\end{split}
\eex
By $W^{2,2}$-estimate of the Lam$\acute{\mbox{e}}$ equation under the Navier-slip boundary
condition (\ref{clcboundary3}) (see \cite{HLX} Lemma 2.3 and also the proof of Lemma 2.2), we obtain,
by using the equation (\ref{clc-2}) and (\ref{upper_bound_p}), 
\beq\label{H2-est}
\begin{split}
\|\nabla^2u\|_{L^{2}}^2\les&\|\mathcal Lu\|_{L^{2}}^2+\|\nabla u\|_{L^2}^2
\les\|\nabla u\|_{L^2}^2+\|\rho\dot{u}\|_{L^{2}}^2+\|\nabla (P(\rho))\|_{L^{2}}^2+\|\de d\cdot\nabla d\|_{L^{2}}^2\\
\les&\|\nabla u\|_{L^2}^2+\|\rho\dot{u}\|_{L^{2}}^2+\|\nabla \rho\|_{L^{2}}^2+\|\nabla d\|_{L^{\infty}}^2\|\de d\|_{L^{2}}^2.
\end{split}
\eeq
Choosing sufficiently small $\epsilon>0$, we have
\bex
\begin{split}
&\frac{d}{dt}\int\left(\mu|\nabla \times u|^{2}+(2\mu+\lambda)|\di u|^{2}+|\nabla\rho|^2+|\de d|^2\right)\,dx+\int\left(2\rho|\dot{u}|^2+|\nabla d_t|^2\right)\,dx\\
\leq&2\frac{d}{dt}\int(M(d):\nabla u-P(\rho)\di u)\,dx\\
+&
C\left[1+\|\mathcal{D} (u)\|_{L^{\infty}}+\|\nabla u\|_{L^{2}}^2+\|\nabla d\|_{L^{\infty}}^2\right]
\left[\|\nabla u\|_{L^{2}}^2+\|\nabla\rho\|_{L^2}^2+\|\de d\|_{L^2}^2\right]+C(1+\|\nabla d\|_{L^\infty}^2).
\end{split}
\eex
Integrating from $0$ to $t$, $0<t<T_*$ and applying (\ref{div-curl}) and (\ref{upper_bound_p}),  
we obtain
\beq\label{ucblp2.40}
\begin{split}
&\int\left(|\nabla u|^{2}+|\nabla\rho|^2+|\de d|^2\right)(t)\,dx
+\int_0^t\int_{\om}\left(\rho|\dot{u}|^2+|\nabla d_t|^2\right)\,dxds\\
\les&\int (|M(d)||\nabla u|+|P(\rho)||\di u|)(t)\,dx+\int (|M(d_0)||\nabla u_0|+|P(\rho_0)||\di u_0|)\,dx\\
+&\int\left(|\nabla u_0|^{2}+|\nabla\rho_0|^2+|\de d_0|^2\right)\,dx+C\int_0^t (1+\|\nabla d\|_{L^\infty}^2)\,ds\\
+&\int_0^t\left[1+\|\mathcal{D} (u)\|_{L^{\infty}}+\|\nabla u\|_{L^{2}}^2+\|\nabla d\|_{L^{\infty}}^2\right]
\left[\|\nabla u\|_{L^2}^2+\|\nabla\rho\|_{L^2}^2+\|\Delta d\|_{L^2}^2\right]\,ds\\
\leq&C+\frac12 \|\nabla u\|_{L^2}^2+C\int (|\nabla d|^4+|\rho|)\,dx\\
+&C\int_0^t\left[1+\|\mathcal{D} (u)\|_{L^{\infty}}+\|\nabla u\|_{L^{2}}^2+\|\nabla d\|_{L^{\infty}}^2\right]
\left[\|\nabla u\|_{L^2}^2+\|\nabla\rho\|_{L^2}^2+\|\Delta d\|_{L^2}^2\right]\,ds.
\end{split}
\eeq
Since the coefficient function
$$\left[1+\|\mathcal{D} (u)\|_{L^{\infty}}+\|\nabla u\|_{L^{2}}^2+\|\nabla d\|_{L^{\infty}}^2\right] 
\in L^1([0,T_*]),$$
the Gronwall's inequality, Lemma 2.3 and the conservation of mass imply that for any $0\le t<T_*$,
$$\int (|\nabla u|^2+|\nabla\rho|^2+|\nabla^2 d|^2)(t)\,dx
+\int_0^t\int_\Omega (\rho |\dot u|^2+|\nabla d_t|^2)\,dx\,ds\le C.$$
The proof is complete.
\endpf\\

As an immediate consequence of the proof of Lemma \ref{uclemma2.4}, we have
\begin{corollary}\label{uccoro2.5}{\it Under the same assumptions 
as in Lemma \ref{uclemma2.4},  we have 
\beq\label{ucblp2.5}
\begin{split}
\sup\ls_{0\leq t<T_*}\left\|d_t\right\|_{L^{2}}^2
+\int_0^{T_*}\left(\|\rho_t\|_{L^2}^2+\|\nabla u\|_{H^1}^2+\|\nabla d\|_{H^2}^2\right)\,dt\leq C.
\end{split}
\eeq}
\end{corollary}

\pf It follows from Lemma \ref{uclemma2.4} that $\nabla u, \nabla\rho, \de d\in L^\infty_tL^2_x(\om\times [0,T_*])$.
By Sobolev's inequality, we then have $u\in L^\infty_tL^6_x(\om\times [0,T_*])$. On the other hand, Lemma \ref{uclemma2.3} implies $\nabla d\in L^\infty_t L^r_x(\om\times [0,T_*])$ for $r=3,4$.  Therefore, by
(\ref{clc-3}), we have
$$|d_t|\les (|u||\nabla d|+|\de d|+|\nabla d|^2)\in L^\infty_tL^2_x(\om\times [0,T_*]).$$
It is easy to see that $L^2_tL^2_x$-estimate of $\nabla^2 u$ and $\nabla^3 d$ follows from 
(\ref{H3-est}) , (\ref{H2-est}), (2.1),  and Lemma \ref{uclemma2.4}. To see $L^2_tL^2_x$-estimate
of $\rho_t$, note that (\ref{clc-1}) and Lemma \ref{uclemma2.2}  imply
$$|\rho_t|\leq |\nabla\rho||u|+\rho|\di u|\les |u||\nabla\rho|+|\nabla u|.$$
By the Sobolev's embedding, we have $u\in L^2_tL^\infty_x(\om\times [0,T_*])$. Hence
$$\||u||\nabla\rho|\|_{L^2(\om\times [0,T_*])}
\le \|u\|_{L^2_tL^\infty_x(\om\times [0,T_*])}\|\nabla\rho\|_{L^\infty_tL^2_x(\om\times [0,T_*])}\le C.$$
This clearly implies $\|\rho_t\|_{L^2(\om\times [0,T_*])}\le C.$ The proof is complete. \endpf\\

\noindent{\bf Step 5}. Estimates of $(\sqrt{\rho}u_t, \nabla^2 u, \nabla d_t, \nabla^3 d)$ in
$L^\infty_tL^2_x(\om\times [0,T_*])$. More precisely, we have 

\blm\label{clclemma2.6} Let $0<T_*<+\infty$ be the maximum time for a strong
solution $(\rho,u,d)$ to (\ref{clc-1})-(\ref{clc-3}), (\ref{clcinitial}) together with (\ref{clcboundary1})
or (\ref{clcboundary2}) or (\ref{clcboundary3}). If (\ref{first3.1}) and (\ref{2.1}) hold, then 
\beq\label{ucblp3.1}
\begin{split}
\sup\ls_{0\leq t<T_*}\int_\om (\rho|u_t|^2+|\nabla^2 u|^2+|\nabla d_t|^2+|\nabla^3 d|^2)\,dx+\int_0^{T_*}\int_{\om}(|\nabla
u_t|^2+|d_{tt}|^2)dxdt\leq C.
\end{split}
\eeq
\elm

\pf For simplicity, we only consider the Navier-slip boundary condition (\ref{clcboundary3}).
Differentiating the equation (\ref{clc-2}) with respect to $t$, we get
\beq\label{ucblp3.2}
\begin{split}
&\rho u_{tt}+\rho_tu_t+\rho u\cdot\nabla u_t+\rho u_t\cdot\nabla
u+\rho_t u\cdot\nabla u+\nabla ((P(\rho))_t)\\
&=\mathcal Lu_t-\nabla\cdot\left(\nabla
d_t\otimes\nabla d +\nabla d\otimes\nabla d_t-\nabla d\cdot\nabla d_t
\mathbb I_3\right).
\end{split}
\eeq 
Since $u_t\cdot\nu=0$ and $\curl u_t\times\nu=0$ on $\pa\om$,
as in the proof of Lemma \ref{uclemma2.4} we can verify
$$-\int\langle\mathcal L u_t, u_t\rangle\,dx=\int (\mu |\nabla\times u_t|^2+(2\mu+\lambda)|\di u_t|^2)\,dx.$$
Thus, multiplying (\ref{ucblp3.2}) by $u_t$ and integrating the resulting equation
over $\om$ and using (\ref{clc-1}), we obtain, by using
Sobolev's inequality, H\"older's inequality, and (\ref{upper_bound_p_t}),
\bex
\begin{split}
&\frac{1}{2}\frac{d}{dt}\int\rho|u_t|^2dx+\int\left(\mu|\nabla \times u_t|^2+(2\mu+\lambda)|\di u_t|^2\right)dx\\
\lesssim&\int\left(\rho |u||\nabla u_t||u_t|+\rho |u||\nabla(u\cdot\nabla
u\cdot u_t)|+|(P(\rho))_t||\di u_t|\right)\,dx\\
&+\int\rho|u_t|^2|\nabla u|\,dx+\int|\nabla d_t||\nabla
d||\nabla u_t|\,dx\\
\lesssim&\|\nabla
u_t\|_{L^2}\|\sqrt{\rho}u_t\|_{L^2}\|u\|_{L^\infty}+\|\nabla
u\|_{L^3}\|u_t\|_{L^6}\|\sqrt{\rho}u_t\|_{L^2}+\|\nabla u_t\|_{L^2}\|\nabla d_t\|_{L^2}\|\nabla
d\|_{L^\infty}\\
&+\|\rho_t\|_{L^2}\|\di
u_t\|_{L^2}+\int(\rho|u||\nabla u|^2|u_t|+\rho |u|^2|\nabla^2 u||u_t|+
\rho|u|^2|\nabla u||\nabla u_t|)\,dx\\
 \lesssim&\|\nabla
u_t\|_{L^2}\|\sqrt{\rho}u_t\|_{L^2}\|\nabla u\|_{H^1}+\|\nabla
u\|_{L^6}\|\nabla u_t\|_{L^2}\|u\|_{L^6}^2+\|u_t\|_{L^6}\|\nabla^2
u\|_{L^2}\|u\|_{L^6}^2\\
&+\|u_t\|_{L^6}\|u\|_{L^6}\|\nabla
u\|_{L^3}^2+\|\rho_t\|_{L^2}\|\di u_t\|_{L^2}+\|\nabla
u_t\|_{L^2}\|\nabla d_t\|_{L^2}\|\nabla d\|_{L^\infty}\\
\lesssim&\|\nabla u_t\|_{L^2}\left(\|\sqrt{\rho}u_t\|_{L^2}\|\nabla
u\|_{H^1}+\|\nabla u\|_{H^1}+\|\rho_t\|_{L^2}+\|\nabla
d_t\|_{L^2}\|\nabla d\|_{L^\infty}\right)\\
\lesssim&\frac{1}{2}\int \mu|\nabla
u_t|^2dx+\|\sqrt{\rho}u_t\|_{L^2}^2\|\nabla
u\|_{H^1}^2+\|\nabla u\|_{H^1}^2+\|\rho_t\|_{L^2}^2+\|\nabla
d_t\|_{L^2}^2\|\nabla d\|_{L^\infty}^2.
\end{split}
\eex This gives
 \beq\label{L^2-estimates of rho u_t}
\begin{split}
&\frac{d}{dt}\int\rho|u_t|^2\,dx+\int\left(\mu|\nabla \times u_t|^2+(2\mu+\lambda)|\di u_t|^2\right)\,dx\\
\lesssim&\|\nabla u\|_{H^1}^2\int\rho|u_t|^2\,dx+\|\nabla^2
u\|_{L^2}^2+\|\rho_t\|_{L^2}^2+\|\nabla
d\|_{L^\infty}^2\|\nabla d_t\|_{L^2}^2+1.
\end{split}
\eeq
Differentiating the equation (\ref{clc-3}) with respect to $t$, multiplying $d_{tt}$ and integrating over $\om$, we obtain,
by using $\frac{\pa d_t}{\pa\nu}=0$ on $\pa\om$, Sobolev and H\"older inequalities, Lemma
\ref{uclemma2.3}, and Lemma \ref{uclemma2.4},
\bex\begin{split}
\frac{1}{2}\frac{d}{dt}\int|\nabla d_t|^2\,dx+\int|d_{tt}|^2\,dx
=&\int\langle\pa_t\left(|\nabla d|^2d-u\cdot\nabla d\right),d_{tt}\rangle \,dx\\
\les&\|d_{tt}\|_{L^2}\|u_t\|_{L^6}\|\nabla d\|_{L^3}+\|d_{tt}\|_{L^2}\|u\|_{L^6}\|\nabla d_t\|_{L^3}\\
&+\|d_{tt}\|_{L^2}\|d_t\|_{L^6}\|\nabla d\|_{L^6}^2+\|d_{tt}\|_{L^2}\|\nabla d_t\|_{L^2}\|\nabla d\|_{L^{\infty}}\\
\leq&\frac{1}{4}\|d_{tt}\|_{L^2}^2+C[\|\nabla u_t\|_{L^2}^2
+\|\nabla d_t\|_{L^2}\|\nabla^2 d_t\|_{L^2}\\
&+(1+\|\nabla d\|_{L^{\infty}}^2)\|\nabla d_t\|_{L^2}^2],
\end{split}
\eex
which implies
\beq\label{ucblp3.6}
\begin{split}
&\frac{d}{dt}\int|\nabla
d_t|^2dx+\int|d_{tt}|^2dx\\ \leq& C[\|\nabla
u_t\|_{L^2}^2+\|\nabla^2 d_t\|_{L^2}\|\nabla d_t\|_{L^2}+(1+\|\nabla d\|_{L^{\infty}}^2)\|\nabla d_t\|_{L^2}^2].
\end{split}
\eeq

Now we need to estimate $\|\nabla^2 d_t\|_{L^2}$. In fact, by applying the standard $H^2$-estimate
on the equation (\ref{clc-3}) and Lemma \ref{uclemma2.3}, we have 
\bex\begin{split}
\|\nabla^2 d_t\|_{L^2}\les&\|\nabla d_t\|_{L^2}+
\| d_{tt}\|_{L^2}+\|\pa_t(u\cdot\nabla d)\|_{L^2}+\|\pa_t(|\nabla d|^2d)\|_{L^2}\\
\les&\|\nabla d_t\|_{L^2}+\| d_{tt}\|_{L^2}+\|u_t\|_{L^6}\|\nabla d\|_{L^3}+\|u\|_{L^6}\|\nabla d_t\|_{L^3}\\
&+\|d_t\|_{L^6}\|\nabla d\|_{L^6}^{2}+\|\nabla d_t\|_{L^3}\|\nabla d\|_{L^6}\\
\les&\| d_{tt}\|_{L^2}+\|\nabla u_t\|_{L^2}
+\|\nabla d_t\|_{L^2}^{\frac{1}{2}}\|\nabla^2 d_t\|_{L^2}^{\frac{1}{2}}+\|\nabla d_t\|_{L^2}\\
\leq&\frac{1}{2}\|\nabla^2 d_t\|_{L^2}+C\left[\| d_{tt}\|_{L^2}+\|\nabla
u_t\|_{L^2}+\|\nabla d_t\|_{L^2}\right].
\end{split}
\eex Thus 
\beq \label{D2t}
\begin{split} \|\nabla^2 d_t\|_{L^2} \les&\|
d_{tt}\|_{L^2}+\|\nabla u_t\|_{L^2}+\|\nabla d_t\|_{L^2}.
\end{split}
\eeq
Substituting (\ref{D2t}) into (\ref{ucblp3.6}), and using
Cauchy inequality, we obtain \bex
\begin{split}
\frac{d}{dt}\int|\nabla
d_t|^2\,dx+\int|d_{tt}|^2\,dx\leq&\frac{1}{4}\|
d_{tt}\|_{L^2}^2
+C\left[\|\nabla
u_t\|_{L^2}^2+(1+\|\nabla
d\|_{L^{\infty}}^2)\|\nabla d_t\|_{L^2}^2\right].
\end{split}
\eex
Thus
\beq\label{ucblp3.7}
\begin{split}
\frac{d}{dt}\int|\nabla
d_t|^2dx+\int|d_{tt}|^2dx\leq&C\left[\|\nabla u_t\|_{L^2}^2+(1+\|\nabla d\|_{L^{\infty}}^2)\|\nabla d_t\|_{L^2}^2\right].
\end{split}
\eeq Multiplying (\ref{ucblp3.7}) by $\frac{\mu}{2C}$ and adding the resulting inequality into
(\ref{L^2-estimates of rho u_t}), applying Lemma \ref{uclemma2.3} and Lemma
\ref{uclemma2.4}, and then employing Gronwall's inequality and applying (\ref{div-curl}) (with $u$ replaced by $u_t$), we obtain
$$\sup_{0\le t<T_*} \int (\rho|u_t|^2+|\nabla d_t|^2)\,dx
+\int_0^{T_*}\int_\Omega(|\nabla u_t|^2+|d_{tt}|^2)\,dx\,dt\le C.$$
To estimate $\nabla^3 d$ in $L^\infty_t L^2_x(\Omega\times [0,T_*])$, first observe that by
Nirenberg's interpolation inequality, we have
$$\|\nabla d\|_{L^\infty}\les \|\nabla d\|_{L^2}+\|\nabla d\|_{L^2}^\frac14\|\nabla^3 d\|_{L^2}^\frac34.$$
Putting this inequality into (\ref{H3-est}) and using $L^\infty_t L^2_x$-bounds of
$\nabla d_t, \nabla u, \nabla^2 d$, we obtain that for any $0\le t<T_*$,
$$\|\nabla^3 d\|_{L^2}^2\leq C+C\|\nabla^3 d\|_{L^2}^\frac32
\leq \frac12\|\nabla^3 d\|_{L^2}^2+ C, $$
which clearly yields that 
$$\sup_{0\le t<T_*}(\|\nabla d\|_{L^\infty}+\|\nabla^3 d\|_{L^2})\le C.$$

To see $\nabla^2 u\in L^\infty_t L^2_x(\om\times [0,T_*])$, observe that
the $H^2$-estimate on the equation (\ref{clc-2}) under (\ref{clcboundary3}),
(\ref{upper_bound_p}), and Lemma \ref{uclemma2.4}  imply that for any $0\le t<T_*$,
 \bex
\begin{split}
\|\nabla^2 u\|_{L^2}\lesssim&\|\nabla u\|_{L^2}+\|\mathcal L u\|_{L^2}\\
\les& 1+\|\sqrt{\rho}u_t\|_{L^2}+\|u\cdot\nabla
u\|_{L^2}+\|\nabla (P(\rho))\|_{L^2}+\|\nabla^2 d\|_{L^2}\|\nabla
d\|_{L^\infty}\\
\lesssim&1+\|\sqrt{\rho}u_t\|_{L^2}+\|u\|_{L^6}\|\nabla
u\|_{L^3}
\lesssim1+\|\sqrt{\rho}u_t\|_{L^2}+\|\nabla
u\|_{L^2}^\frac{1}{2}\|\nabla^2 u\|_{L^2}^\frac{1}{2}\\
\leq&\frac{1}{2}\|\nabla^2u\|_{L^2}^2+C.
\end{split}
\eex 
In particular, we have
$$\sup_{0\le t<T_*}\|\nabla^2 u\|_{L^2}\le C.$$
The proof is now complete.
\endpf\\

\noindent{\bf Step 6}. Estimate of $\nabla \rho$ in $L^\infty_t L^q_x(\om\times [0,T_*])$ for
some $3<q\le 6$. With the estimates already established by the previous Lemmas, we then
have the following Lemma.
\blm\label{clclemma2.7} Let $0<T_*<+\infty$ be the maximum time for a strong
solution $(\rho,u,d)$ to (\ref{clc-1})-(\ref{clc-3}), (\ref{clcinitial}) together with (\ref{clcboundary1})
or (\ref{clcboundary2}) or (\ref{clcboundary3}). If (\ref{first3.1}) and (\ref{2.1}) hold, then 
\beq\label{ucblp3.8}
\begin{split}
\sup\ls_{0\leq t<T_*}\left(\max_{r=2, q}\|\rho_t\|_{L^r}+\|\rho\|_{W^{1,q}}\right)+\int_0^{T_*}\big(\|u\|_{D^{2,q}}^2+\|\nabla^2d_t\|_{L^2}^2+\|\nabla^4
d\|_{L^2}^2\big)\,dt\leq C,
\end{split}
\eeq for any $3<q\leq 6$.  \elm

\pf For $3<q\le 6$,  by the same calculations as in \cite{HLX} Lemma 2.5, we have
\bex
(|\nabla\rho|^q)_t+\di(|\nabla\rho|^qu)+(q-1)|\nabla\rho|^q\di
u+q|\nabla\rho|^{q-2}(\nabla\rho)^t\mathcal D(u)(\nabla\rho)+q\rho|\nabla\rho|^{q-2}\nabla\rho\cdot\nabla\di
u=0, \eex
which yields that for (\ref{clcboundary1}) or (\ref{clcboundary2}) 
\beq\label{grad 0f
density-1} \frac{d}{dt}\|\nabla\rho\|_{L^q}\le
C\left(\|\mathcal D(u)\|_{L^\infty}+1\right)\|\nabla\rho\|_{L^q}+C\|\nabla\di
u\|_{L^q}, \eeq and that for (\ref{clcboundary3})
 \beq\label{grad 0f density-2}
\frac{d}{dt}\|\nabla\rho\|_{L^q}\le
C\left(\|\mathcal D(u)\|_{L^\infty}+1\right)\|\nabla\rho\|_{L^q}+C\|\nabla
G\|_{L^q}, \eeq where $G\equiv (2\mu+\lambda)\di u-P(\rho)$.

For boundary conditions (\ref{clcboundary1}) or (\ref{clcboundary2}),  by using the $L^p$-estimate for the
elliptic equation and (\ref{upper_bound_p}) we have
\beq\label{W2,q of u}
\begin{split}
\|\nabla^2u\|_{L^q}\lesssim& \|\rho u_t\|_{L^q}+\|\rho u\cdot\nabla
u\|_{L^q}+\|\nabla (P(\rho))\|_{L^q}+\|\Delta d\cdot\nabla d\|_{L^q}+1\\
\lesssim&\|\sqrt{\rho}u_t\|_{L^2}^\frac{6-q}{2q}\|u_t\|_{L^6}^\frac{3q-6}{2q}+\|u\|_{L^\infty}\|\nabla
u\|_{L^q}+\|\nabla\rho\|_{L^q}+\|\nabla d\|_{L^\infty}\|\Delta
d\|_{L^q}+1\\ \lesssim& \|\nabla u_t\|_{L^2}+\|\nabla\rho\|_{L^q}+1.
\end{split}
\eeq Substituting (\ref{W2,q of u}) into (\ref{grad 0f density-1}),
and using Gronwall's inequality, we obtain the bound $\sup\limits_{0\le t<T_*}\|\rho\|_{W^{1,q}}$
for the first two boundary conditions (\ref{clcboundary1}) and (\ref{clcboundary2}).

For boundary condition (\ref{clcboundary3}), we rewrite (\ref{clc-2}) 
as \beq\label{Grad G}
\nabla G=\mu\nabla\times\mathrm{curl}u+\rho u_t+\rho u\cdot\nabla
u+\Delta d\cdot\nabla d, \eeq 
which yields that $G$
satisfies \beq\label{de G} \begin{cases}\de G=\di \left(\rho
u_t+\rho u\cdot\nabla u+\nabla d\cdot\Delta d\right),\mbox{ in} \ \Omega,\\
\nabla G\cdot\nu=-\rho(u\cdot\nabla)\nu\cdot u,\mbox{ on }\ \partial\Omega,
\end{cases} \eeq
where we have used that
$(\nabla\times\mathrm{curl}u)\cdot\nu|_{\partial\Omega}=0$
\big($\mathrm{curl}u\times\nu|_{\partial\Omega}=0\ \mathrm{implies}\
(\nabla\times\mathrm{curl}u)\cdot\nu|_{\partial\Omega}=0$, see \cite{HLX} page 33 or \cite{chen-osborne-qian1,
chen-osborne-qian2}), $\nabla
d\cdot\nu|_{\partial\Omega}=0$ and $u\cdot\nu|_{\partial\Omega}=0$.

Using the $L^p$-estimate for Neumann problem to the elliptic equation (\ref{de G}), we have
\beq\label{LpgradG}
\begin{split} \|\nabla G\|_{L^q}\lesssim&\|\rho u_t\|_{L^q}+\|\rho
u\cdot\nabla u\|_{L^q}+\|\nabla d\cdot\Delta
d\|_{L^q}+\|\rho|u|^2\|_{C(\overline{\Omega})}\\
\lesssim&\|\nabla u_t\|_{L^2}+1.
\end{split} \eeq
Puting (\ref{LpgradG}) into (\ref{grad 0f density-2}),  we  obtain the bound $\sup\limits_{0\le t<T_*}\|\rho\|_{W^{1,q}}$
by Gronwall's inequality.

For $r=2$ or $q$, (\ref{clc-1}) implies \bex
\|\rho_t\|_{L^r}&\lesssim&
\|u\|_{L^\infty}\|\nabla\rho\|_{L^r}+\|\rho\|_{L^\infty}\|\di
u\|_{L^r}\\ &\lesssim&\|\nabla
u\|_{H^1}\|\nabla\rho\|_{L^r}+\|\rho\|_{L^\infty}\|\nabla
u\|_{H^1}\le C. \eex
It follows from (\ref{D2t}) that
\bex \|\nabla^2d_t\|_{L^2}^2\lesssim \|d_{tt}\|_{L^2}^2+\|\nabla
u_t\|_{L^2}^2+1. \eex 
This, with the help of (\ref{ucblp3.7}), implies, after integrating over $[0,T_*]$, 
\beq\label{grad^2d_t} \int_0^T\|\nabla^2d_t\|_{L^2}^2dt\le C.
\eeq
Applying the standard $L^2$-estimate to  (\ref{clc-3}), we have
\bex
\|\nabla^4d\|_{L^2}^2&\lesssim&\|\nabla^2d_t\|_{L^2}^2+\|\nabla^2(u\cdot\nabla
d)\|_{L^2}^2+\|\nabla^2(|\nabla d|^2d)\|_{L^2}^2+1\\
&\lesssim&\|\nabla^2d_t\|_{L^2}^2+\|u\|_{L^\infty}^2\|\nabla^3d\|_{L^2}^2+\|\nabla
d\|_{L^\infty}^2\|\nabla^2u\|_{L^2}^2+\|\nabla
u\|_{L^6}^2\|\nabla^2d\|_{L^3}^2+1\\
&\lesssim&\|\nabla^2d_t\|_{L^2}^2+1. \eex
Integrating this inequality over $[0,T_*]$, and using (\ref{grad^2d_t}), we get \bex
\int_0^T\|\nabla^4d\|_{L^2}^2dt\le C. \eex 
By the bound on $\|\nabla\rho\|_{L^q}$ in (\ref{ucblp3.8}),
(\ref{W2,q of u}) and (\ref{ucblp3.1}), we easily see that  
\bex \int_0^{T_*}\|\nabla^2u\|_{L^q}^2\,dt\le C.
\eex
holds for (\ref{clcboundary1}) or (\ref{clcboundary2}). For the boundary condition (\ref{clcboundary3}), since
$u\cdot\nu=0$ on $\pa\om$, it follows from Bourguignon-Brezis \cite{BB} (see also \cite{HLX} Lemma 2.3) and (\ref{upper_bound_p}) that
\bex
\|\nabla^2 u\|_{L^q} &\les&\|\nabla (\di u)\|_{L^q}+\|\nabla(\curl u)\|_{L^q}+\|\nabla u\|_{L^q}\\
&\les&\|\nabla G\|_{L^q}+\|\nabla\rho\|_{L^q}+\|\nabla u\|_{H^1}+\|\nabla(\curl u)\|_{L^q}\\
&\les& 1+\|\nabla G\|_{L^q}+\|\nabla(\curl u)\|_{L^q}.\\
\eex
Since $(\nabla\times u)^\tau=0$ on $\pa\om$, it follows from \cite{Wahl} that
\bex
\|\nabla(\curl u)\|_{L^q}&\les&\|\di(\curl u)\|_{L^q}+\|\nabla\times \curl u\|_{L^q}
\les
\|\nabla\times \curl u\|_{L^q},
\eex
where we have used the fact that $\di(\curl u)=0$.  On the other hand, since
$$\mu  \nabla\times \curl u=\nabla G-\rho u_t-\rho u\cdot\nabla u-\de d\cdot\nabla d,$$
(\ref{LpgradG}) implies
$$\|\nabla\times \curl u\|_{L^q}\les 1+\|\nabla u_t\|_{L^2}.$$
Putting these estimates together, we have
$$\|\nabla^2 u\|_{L^q}\les 1+\|\nabla u_t\|_{L^2}+\|\nabla G\|_{L^q},$$
which clearly implies
$$\int_0^{T_*}\|\nabla^2 u\|_{L^q}^2\le C.$$
The proof is now complete. 
\endpf\\

\noindent{\bf Step 7}. Completion of proof of Theorem \ref{umaintheorem}: \\

With the above established estimates, we obtain (2.2) and (2.3). This implies that
$T_*$ is not the maximum time of existence of strong solutions, which
contradicts the definition of $T_*$.
Therefore, (\ref{2.1}) is false. The proof of Theorem \ref{umaintheorem} is now complete.
\qed

\end{document}